\documentclass[12pt]{article}
\usepackage{amsfonts}

\newtheorem{thm}{Theorem}[section]
\newtheorem{lemma}[thm]{Lemma}
\newtheorem{cor}[thm]{Corollary}

\newenvironment{remark}{\par\medskip\noindent{\bf Remark.\ }}{\par\smallskip}
\newcommand{\proof
}{\par\medskip\noindent {\bf Proof.\ \ }}

\newcommand{\be}{\begin{equation}}
\newcommand{\ee}{\end{equation}}
\newcommand{\openbox}{\leavevmode
  \hbox to8pt{\hfil\vrule\vbox to6pt{\hrule width6pt\vfil\hrule}\vrule}}

\newcommand{\qed}{\hbox to5pt{ } \hfill \openbox\bigskip\medskip}

\newcommand{\Zp}{\mathbb Z _p}
\newcommand{\Zk}{\mathbb Z _k}

\newcommand{\Zq}{\mathbb Z _q}

\newcommand{\cF}{\mbox{$\cal F$}}

\newcommand{\cP}{\mbox{$\cal P$}}

\newcommand{\N}{\mathbb N}
\newcommand{\Z}{\mathbb Z}

\newcommand{\R}{\mathbb R}

\title{A new exponential upper bound for the Erd\H{o}s-Ginzburg-Ziv constant}
\author{G\'abor Heged\"{u}s
\\{\normalsize }
}

\begin{document}
\maketitle

\begin{abstract}
Naslund used Tao's slice rank bounding method to give new exponential upper bounds for the   Erd\H{o}s--Ginzburg-Ziv  constant of finite Abelian groups of high rank. In our short manuscript we improve slightly Naslund's upper bounds. We extend Naslund's results and prove new exponential upper bounds for the   Erd\H{o}s--Ginzburg-Ziv  constant of arbitrary finite Abelian groups. Our main results depend on a conjecture about Property D. 
\end{abstract}
\medskip

\noindent

\section{Introduction}

%3.  Sets $A\subseteq G$ satisfying the condition 
%$$
%0\notin hA
%$$
%for each $h\geq 1$ are called {\em zero-sum-free} sets. 

Let $A$ denote an additive finite Abelian group.
Let $\exp(A)$ denote the exponent of $A$.

%Let $D(A)$ denote the smallest integer $\ell\in \N$ such that every sequence $S$ over $G$ of length $|S|\geq \ell$ has a non-trivial zero--sum subsequence.

%Then $D(A)$ is the {\em Davenport} constant of $A$.

We denote by $\eta(A)$ the smallest integer $\ell\in \N$ such that such that every sequence $S$ over $G$ of length $|S|\geq \ell$ has a zero--sum sub-sequence of length $1\leq |T|\leq exp(A)$.

We denote by $s(A)$ the smallest integer $\ell\in \N$ such that every sequence $S$ over $G$ of length $|S|\geq \ell$ has a zero--sum sub-sequence of length $|T|=exp(A)$.

Then $s(A)$ is the {\em Erd\H{o}s-Ginzburg-Ziv} constant of $A$.

%Let $g(A)$ denote the smallest integer $\ell\in \N$ such that every square-free sequence $S$ over $G$ of length $|S|\geq \ell$ has a zero--sum sub-sequence of length $|T|=exp(A)$.

%The precise value of $s(A)$  is known for groups with rank at most two. We have
%\begin{thm}
%If $A=\Z_{n_1} \oplus \Z_{n_2}$, where $1\leq n_1 | n_2$, then
%$s(A)=2n_1+2n_2-3$. 
%\end{thm}
We use frequently the following result (see \cite{CDGGS} Proposition 3.1).

\begin{thm} \label{exp_upper}
Let $G$ be a finite Abelian group and let $H\leq G$ be a subgroup such that $exp(G)=exp(H)exp(G/H)$. Then
$$
s(G)\leq exp(G/H)(s(H)-1)+s(G/H).
$$
\end{thm}

The following Lemma will be useful in our proofs (see \cite{FS} Lemma 3.5).
\begin{lemma} \label{Foxlem}
Let $A$ be a finite Abelian group. Let us write $A$ as
$$
A \cong A(p_1) \times \ldots \times A(p_m)
$$
where  each $A_i:=A(p_i)$ is a $p_i$-group. Then  each $A_i$ can be written as  a 
product of cyclic groups whose orders are power of $p_i$. Let $n_i$ denote the number of these cyclic factors. Then
$$
s(A)< exp(A) \Big( \sum_{j=1}^m \frac{s({\Z}_{p_j}^{n_j})}{p_j-1}\Big).
$$

\end{lemma}
%Class res:

The following inequality is well-known, see \cite{H}.
\begin{thm} \label{Har}
Let $k\geq 2$, $n\geq 1$ be arbitrary integers.
Let $A:={(\Zk})^n$. Then
$$
(k-1) 2^n +1\leq s(A)\leq (k-1)k^n+1.
$$
\end{thm}

%The following result was proved in \cite{EEGKR} as Lemma 2.3(3).
%\begin{thm} \label{PropD2}
%Let $A$ denote a finite Abelian group with $exp(A)=k\geq 2$. Then 
%$$
%g(A)\leq s(A)\leq (g(A)-1)(k-1)+1.
%$$
%If $A={(\Zk})^n$ with $k\geq 2$, $n\in \N$ and $s(A)= (g(A)-1)(k-1)+1$, then $A$ has Property D.
%\end{thm}

Harborth determined $s(A)$ in the following special case  in  \cite{H}.
\begin{thm} \label{Har2}
Let $a\geq 1$, $n\geq 1$ be arbitrary integers. Let $k:=2^a$ and $A:={(\Zk})^n$. Then 
$$
s(A)=(k-1) 2^n +1.
$$
\end{thm}

Let $A:=({\Zk})^n$ with $k,n\in \N$ and $k\geq 2$. We can ask for the structure of sequences of length $s(A)-1$ that do not have a zero-sum sub-sequence of length $k$. 
The following conjecture is well--known:  every group  $A:=({\Zk})^n$ satisfies  Property D (see \cite{GG}, Conj. 7.2).

Property D: Every sequence $S$ over $A$ of length $|S|=s(A)-1$ that has no zero-sum sub-sequence of length $k$ has the form $S=T^{k-1}$ for some subset $T$ over $A$.

In the following Theorem we collected all known groups satisfying  Property D. 
\begin{thm} \label{PropD}
The following groups has Property D:
\begin{itemize}
\item[(i)] $A=({\Zk})^n$, where $k=2^{\alpha}$, $\alpha,n\geq 1$ is arbitrary;
\item[(ii)] $A=({\Zk})^n$, where $k=3$, $n\geq 1$ is arbitrary (see \cite{H}, Hilfsatz 3);
\item[(iii)] $A=({\Zk})^n$, where $n=1$, $k\geq 2$ is arbitrary;
\item[(iv)] $A=({\Zk})^n$, where $n=2$, $k$  is not divisible a prime greater than $7$ (see \cite{ST});
\item[(v)] $A=({\Zk})^n$, where $n=3$, $k=5^{\alpha}$, ${\alpha}>0$ (see \cite{GHST}, Theorem 1.9);
\item[(vi)] $A=({\Zk})^n$, where $n=3$, $k=3^{\alpha}$, $\alpha>0$ (see \cite{GT} Corollary 1.1).
\end{itemize}
\end{thm}

%Lower bounds:
Elsholtz proved the following lower  bound for $s(A)$ in \cite{E}, where 
$A:=({\Zk})^n$.
\begin{thm} \label{Els}
Let $k$ be an odd integer.  The following inequality holds:
$$
s(({\Zk})^n)\geq (1.125)^{\lfloor \frac{n}{3}\rfloor}(k-1)2^n+1.
$$
\end{thm}

These remarkable lower bounds appeared in \cite{EEGKR}:
\begin{thm} \label{low}
Let $k$ be an odd integer. Then 
$\eta(({\Zk})^3)\geq 8k-7$ and $s(({\Zk})^3)\geq 9k-8$.
\end{thm}

\begin{thm} \label{low2}
Let $k$ be an odd integer with $k\geq 3$. Then 
$\eta(({\Zk})^4)\geq 19k-18$ and $s(({\Zk})^4)\geq 20k-19$.
\end{thm}

Let $A:=({\Zk})^r$. Let $\cP$ denote the set of all prime factors of $k$. One of our main result is a better bound for $s(A)$ if we suppose that Property D is satisfied for all groups $({\Zp})^n$, where $p\in \cP$. We give also new exponential upper   bounds for the numbers $s(({\Zq})^n)$, where $q$ is an arbitrary prime  power.

%Chintamani, Moriya, Gao,  Paul and Thangadurai proved in \cite{CMGPT} the following upper bounds for the Erd\H{o}s-Ginzburg-Ziv and the Davenport constants.

%\begin{thm} \label{CMGPT1}
%Let $A$ be an arbitrary finite abelian group of rank $r$ with invariants $n_1, n_2, \ldots n_r$. 
%Then
%$$
%D(A)\leq n_r+n_{r-1}+(c(3)-1)n_{r-2}+(c(4)-1)n_{r-3}+\ldots + (c(r)-1)n_1 +1.
%$$
%\end{thm}

%\begin{thm} \label{CMGPT2}
%Let $A$ be an arbitrary finite abelian group of rank $r$ with invariants $n_1, n_2, \ldots n_r$.
%Then
%$$
%s(A)\leq n_rc(1)+n_{r-1}c(2)+c(3)n_{r-2}+c(4)n_{r-3}+\ldots + c(r)n_1 +1.
%$$
%\end{thm}

%\begin{thm} \label{CMGPT3}
%Let $A$ be an arbitrary finite abelian group of rank $r$ with invariants $n_1, n_2, \ldots n_r$. 
%Then
%$$
%\eta(A)\leq n_r(c(1)-1)+n_{r-1}(c(2)-1)+(c(3)-1)n_{r-2}+(c(4)-1)n_{r-3}+\ldots + (c(r)-1)n_1 +1.
%$$

%\end{thm}

Naslund achieved the following breakthrough in \cite{N} Theorem 2. 
\begin{thm}
Let $k\geq 2$ be a fixed integer. Let  $q$ denote the largest prime power dividing $k$.
Suppose that $A:=({\Zk})^n$ satisfies Property D. Then
$$
s(A)\leq (k-1)(\gamma_{k,q})^n +1,
$$
where
$$
\gamma_{k,q}=\frac{k}{q}\inf_{0<x<1} \frac{1-x^q}{1-x} x^{-\frac{q-1}{k}}.
$$
In particular, if $q$ is a  prime power and  $A:=({\Zq})^n$ satisfies Property D, then 
$$
s(A)\leq (q-1)4^n+1.
$$
\end{thm}

\begin{remark}
The manuscript \cite{N} contains some typos in Theorem 2. 
\end{remark}
\begin{remark}
It is easy to check that
$$
\gamma_{q,q}\sim \Big( \frac{2^q-1}{2^q}\Big)2^{\frac{2q-1}{q}},
$$
when $q$ is a sufficiently large prime power.
\end{remark}
%The following result was proved in \cite{EEGKR} as Theorem 1.4.

%\begin{thm} \label{EGZupper}
%Let $A:={\Z}_{n_1}\bigoplus \ldots \bigoplus {\Z}_{n_r}$, where $r:=r(G)$ and $1<n_1|\ldots | n_r$. Let $c_1, \ldots ,c_r\in \N$ such that for all primes $p\in \cP$ with $p| n_r$ and all $1\leq i\leq r$ we have 
%$$
%s({\Zp}^i)\leq c_i(p-1)+1.
%$$
%Then
%$$
%s(A)\leq \sum_{i=1}^r (c_{r+1-i}-c_{r-i})n_i -c_r +1,
%$$
%where $c_0=0$. In particular, if $n_1=\ldots =n_r=n$, then $s(A)\leq c_r(n-1)+1$.
%\end{thm}

%The proof is based on Theorem \ref{exp_upper}.

%\proof
Finally we use the following well--known Lemma in the proofs of our main results.
\begin{lemma} \label{monom}
Consider the set of monomials 
$$
B(n,k):=\{x_1^{\alpha_1}\cdot \ldots \cdot x_n^{\alpha_n}:~ \sum_i \alpha_i \leq k\}.
$$
Then
$$
|B(n,k)|={n+k \choose n}.
$$

\end{lemma}
%Consider the prime factorization of $m$ as $m=p_1^{t_1}\ldots p_s^{t_s}$. Now we can prove Corollary \ref{maincor2} by induction on $s$.

%Clearly if $s=1$, then we get the result by Theorem \ref{ppower}. Suppose that our inequality (\ref{genEGZ}) is true for  a  fixed $s$. We will prove that (\ref{genEGZ}) is true for $s+1$. 
%\begin{cor} \label{maincor4}
%Let  $u\geq 3$ be an integer. Let $p_i$ denote the prime factors of $u$, where $1\leq i\leq s$. 
%Suppose that the groups $A_i:=({\Z}_{p_i})^n$ satisfy Property D for each $1\leq i\leq s$.

%Then there exists a constant $1<d(u)<u$ depending only on $u$ such that
%\be 
%s(({\Z}_u)^n)\leq d(u)^n
%\ee
%for each $n\geq 1$.
%\end{cor}

\section{Main results}

We state here our main results.

\begin{thm} \label{main}
Let $p$ be a prime.  Let $n\geq 1$ be an integer. Suppose that Property D is satisfied for the group $({\Zp})^n$.
Then
$$
s(({\Zp})^n)\leq p(p-1){\frac{n(2p-1)}{p}\choose n}.
$$
\end{thm}
\begin{cor} \label{mainc1}
Let $p$ be a fixed prime.  Let $n\geq 1$ be an integer. Suppose that Property D is satisfied for the group $({\Zp})^n$.
Then
$$
s(({\Zp})^n)\leq (p-1)\Big(\Big( 2+ \frac{1}{p-1}\Big)^{\frac{p-1}{p}}\Big(2-\frac{1}{p} \Big)\Big)^n+1.
$$ 
\end{cor}

\begin{cor} \label{mainc2}
Let $p$ be a fixed prime.  Let $n\geq 1$ be an integer. Suppose that Property D is satisfied for the group $({\Zp})^n$.
Then
$$
s(({\Zp})^n)\leq p(p-1){2n \choose n}+1.
$$ 
\end{cor}

\begin{remark}
It is easy to check from Stirling's formula that
$$
{2n \choose n} \sim \frac{4^n}{\sqrt{\pi n}},
$$
when $n$ is sufficiently large.
\end{remark}
\begin{thm} \label{main2}
Let $q=p^{\alpha}\geq 3$ be an odd prime power. Let $n\geq 1$ be an integer. Suppose that Property D is satisfied for the group $({\Zp})^n$.
Then
$$
s(({\Zq})^n)\leq p(q-1){2n \choose n}+1.
$$ 

\end{thm}

We can extend Theorem  \ref{main2} from a prime power to an arbitrary  composite number.
\begin{thm} \label{main3}
Let $k\geq 2$ be a fixed odd integer. We can factorize $k$ as 
$$
k=p_1^{\alpha_1} \cdot\ldots \cdot p_r^{\alpha_r}
$$
where $p_i$ are distinct primes.

Let $\cP$ denote the set of all prime factors of $k$. 
Suppose that Property D is satisfied for each groups $({\Zp})^n$, where $p\in \cP$. Then
$$
s(({\Zk})^n)\leq (p_1\cdot\ldots \cdot p_r)(k-1){2n \choose n}+1.
$$ 

\end{thm}
%\proof 

\begin{thm} \label{main4}
Let $A$ be a finite Abelian group. We can write $A$ as
$$
A \cong A(p_1) \times \ldots \times A(p_m)
$$
where  each $A_i:=A(p_i)$ is a $p_i$-group. Then  each $A_i$ is  a 
product of cyclic groups whose orders are power of $p_i$. Let $n_i$ denote the number of these cyclic factors. Suppose that Property D is satisfied for each  groups $({\Z}_{p_i})^{n_i}$, where $1\leq i\leq m$. Then
$$
s(A)< exp(A) \Big(  \sum_{j=1}^m  p_j{2n_j \choose n_j} + \sum_{j=1}^m \frac{1}{p_j-1}\Big).
$$
\end{thm}

\begin{thm} \label{main5}
Let $k:=3^{\alpha}5^{\beta}$, where $\alpha, \beta\geq 0$, $\alpha+\beta\geq 1$
are integers. Then
$$
s(({\Zk})^3)\leq 300k-299
$$
and
$$
\eta(({\Zk})^3)\leq 299k-298.
$$
\end{thm}

\section{Proofs of the main results}

{\bf Proof of Theorem \ref{main}:}\\

First we prove the following Theorem.
\begin{thm} \label{subset}
Suppose that $A\subseteq ({\Zp})^n$ satisfies
$$
|A| > p {\frac{n(2p-1)}{p}\choose n}.
$$ 

Then $A$ contains $p$ not necessarily distinct but not all equal elements $v_1, \ldots ,v_p$ such that
$$
\sum_i v_i=0.
$$
\end{thm}
\proof 
Indirectly, suppose that $A$ doesn't contain $p$ not necessarily distinct but not all equal elements $v_1, \ldots ,v_p$ such that
$$
\sum_i v_i=0.
$$

Then it follows from Tao's slice rank bounding method (see \cite{T}, \cite{N} Proposition 1 and inequality 4.2) that 
$$
|A|\leq p\cdot |\{x_1^{\alpha_1}\cdot\ldots \cdot x_n^{\alpha_n}:~ 0\leq \alpha_i \leq p-1 \mbox{ for each } i,\ \sum_i \alpha_i \leq \frac{n(p-1)}{p}  \}|.
$$

But
$$
|\{x_1^{\alpha_1}\cdot\ldots\cdot x_n^{\alpha_n}:~ 0\leq \alpha_i \leq p-1 \mbox{ for each } i,\ \sum_i \alpha_i \leq \frac{n(p-1)}{p}  \}|\leq 
$$ 
$$
|\{x_1^{\alpha_1}\cdot\ldots\cdot x_n^{\alpha_n}:~ \sum_i \alpha_i \leq \frac{n(p-1)}{p} \}|={\frac{n(2p-1)}{p}\choose n}
$$
by Lemma \ref{monom}, hence
$$
|A|\leq p {\frac{n(2p-1)}{p}\choose n}.
$$
\qed

Theorem \ref{main} follows easily from the assumption that Property D is satisfied for the group $({\Zp})^n$ and Theorem \ref{subset}. 

Namely let $S$ be an arbitrary sequence in  $({\Zp})^n$ of length $ s( ({\Zp})^n)-1$ for which there exist no $p$ elements that sum to zero. Then Property D implies that we can write $S$ as a multi-set in the form
$$
S=\cup_{i=1}^{p-1} B,
$$
where $B \subseteq ({\Zp})^n$ is a subset. Clearly $B$ doesn't contain $p$ not necessarily distinct  but not all equal elements that sum to zero.

We get from Theorem \ref{subset}  that
$$
|B|\leq p {\frac{n(2p-1)}{p}\choose n},
$$
consequently
$$
s( {\Zp}^n)\leq p (p-1){\frac{n(2p-1)}{p}\choose n}+1.
$$
\qed

{\bf Proof of Corollary \ref{mainc1}:}\\
First we prove:

\begin{thm} \label{subset2}
Suppose that $A\subseteq ({\Zp})^n$ satisfies
$$
|A| > \Big(\Big( 2+ \frac{1}{p-1}\Big)^{\frac{p-1}{p}}\Big(2-\frac{1}{p} \Big)\Big)^n.
$$ 

Then $A$ contains $p$ not necessarily distinct but not all equal elements $v_1, \ldots ,v_p$ such that
$$
\sum_i v_i=0.
$$
\end{thm}
\proof 
Indirectly, suppose that $A$ doesn't contain $p$ not necessarily distinct but not all equal elements $v_1, \ldots ,v_p$ such that
$$
\sum_i v_i=0.
$$
Then we get from the proof  Theorem of  \ref{subset} that
$$
|A|\leq  p{\frac{n(2p-1)}{p}\choose n},
$$
Sondow and Zudilin proved the following simple upper bound for the binomial coefficients in \cite{SZ}:

\begin{thm} \label{Bin_upper}
Let $m\geq 1$ be a positive integer and $r\in \R$ be an arbitrary real. Then
$$
{(r+1)m\choose m}\leq \Big( \frac{(r+1)^{r+1}}{r^r} \Big)^m.
$$
\end{thm}
\qed

It follows from  Theorem \ref{Bin_upper} with the choices 
$m:=n$ and $r:=\frac{n(p-1)}{p}$ that 
$$
|A|\leq p\Big(\Big( 2+ \frac{1}{p-1}\Big)^{\frac{p-1}{p}}\Big(2-\frac{1}{p} \Big)\Big)^n
$$

Finally it follows from a standard amplification argument, that 
\begin{equation} \label{upper2}
|A| \leq  \Big(\Big( 2+ \frac{1}{p-1}\Big)^{\frac{p-1}{p}}\Big(2-\frac{1}{p} \Big)\Big)^n,
\end{equation}
which gives a contradiction. 

Namely let $m\geq 1$ be arbitrary and  consider the set $A^m\subseteq (\Zp)^{nm}$. Then 
$$
|A|^m=|A^m|\leq p\Big(\Big( 2+ \frac{1}{p-1}\Big)^{\frac{p-1}{p}}\Big(2-\frac{1}{p} \Big)\Big)^{nm}.
$$
Consequently
$$
|A| \leq  p^{1/m}\Big(\Big( 2+ \frac{1}{p-1}\Big)^{\frac{p-1}{p}}\Big(2-\frac{1}{p} \Big)\Big)^n,
$$ 
and if $m$ tends to infinity, then we get the inequality (\ref{upper2}).

\qed

Corollary \ref{mainc1} is a clear consequence of the assumption that Property D is satisfied for the group $({\Zp})^n$ and Theorem \ref{subset2}. \qed

{\bf Proof of Corollary \ref{mainc2}:}\\

Corollary \ref{mainc2} follows obviously from Theorem \ref{main}.
\qed

{\bf Proof of Theorem \ref{main2}:}\\
Let $G:=({\Zq})^n$ and $H:=({\Zp})^n$. Clearly
$$
G/H \cong (\Z_{p^{\alpha-1}})^n
$$
and $\exp(G/H)=p^{\alpha-1}$.
%\section{Preliminaries}

By the induction hypothesis we get
$$
s(G/H)\leq p(p^{\alpha-1}-1){2n \choose n}+1.
$$
It follows from Theorem \ref{mainc2} that
$$
s(H)\leq  p(p-1){2n \choose n}+1.
$$

We can apply Theorem \ref{exp_upper} for $G$ and $H$, since $exp(G)=exp(H)exp(G/H)$:
$$
s(G)\leq exp(G/H)(s(H)-1)+s(G/H)\leq 
$$
$$
\leq p^{\alpha-1}\cdot \Big( p(p-1){2n \choose n} \Big)+(p^{\alpha}-p){2n \choose n}+1=p(p^{\alpha}-1){2n \choose n}+1.
$$

\qed

{\bf Proof of Theorem \ref{main3}:}\\

We can derive Theorem \ref{main3}  easily from Theorem \ref{exp_upper}. This proof is very similar to the proof of Theorem \ref{main2}, so we omit it.
\qed

{\bf Proof of Theorem \ref{main4}:}\\

It follows from Theorem \ref{mainc2} that 
$$
s({\Z}_{p_j}^{n_j})\leq p_j(p_j-1){2n_j \choose n_j}+1
$$
for each $1\leq i\leq m$. If we combine this result with Lemma \ref{Foxlem}, then we get our result.

\qed

{\bf Proof of Theorem \ref{main5}:}\\

Theorem \ref{main5} follows clearly from Theorem \ref{PropD} (v) and (vi) and Theorem \ref{main3}.
\qed

{%\bf Acknowledgements.} 

\end{document}